\documentclass[12pt]{article}
\usepackage[utf8]{inputenc}
\usepackage[T1]{fontenc}
\usepackage{amssymb}
\usepackage{amsthm}
\usepackage{amsmath}
\usepackage{graphicx}
\usepackage{color}
\usepackage{amsfonts}
\usepackage{epsfig}
\usepackage{latexsym}
\usepackage{graphicx}
\usepackage{pdflscape}
\usepackage{longtable}
\usepackage{lscape}
\usepackage{bm}
\usepackage{subfigure}
\usepackage{tabulary}
\usepackage{tabularx}
\usepackage{multirow}
\usepackage{graphics}
\usepackage{relsize}
\usepackage{enumitem}
\newtheorem{theorem}{Theorem}[section]
\newtheorem{lemma}[theorem]{Lemma}

\usepackage[parfill]{parskip}
\usepackage{authblk}

\title{A Note on NBUE and NWBUE \\Classes of Life Distributions }
\author{M. Z. Anis\thanks{zafar@isical.ac.in}, }

\affil{SQC \& OR Unit, Indian Statistical Institute, \\ 203 B T Road,Calcutta 700 108, India}

\begin{document}

\maketitle

\begin{abstract}
Non-monotonic ageing notions are looked upon as an extension of the corresponding monotonic ageing notions in this work. In particular,   the New Better than Used in Expectation (NBUE) and the corresponding non-monotonic analogue New Worse then Better than Used in Expectation (NWBUE) classes of life distributions is considered. Some additional results for the NBUE class are obtained. While many properties of the NBUE class carry over in an analogous way to the NWBUE class, it is shown by means of counterexamples that  the moment bounds do not. Some corrective results with respect to popular notions of the NWBUE class are also presented. 
\end{abstract}


\emph{Secondary:} 60E99; 62N05

\section{Introduction}
In the context of reliability theory, the most popularly used distribution is the exponential distribution, primarily due to its elementary  mathematical structure. It essentially
characterizes ``\emph{no ageing}'', which in many cases may not be
correct. The alternatives could involve monotonic ageing or non-monotonic
ageing notions. Monotonic ageing means that the same, i.e. identical type of ageing
continues throughout the entire lifetime of the unit/component. Such ageing could
be positive (meaning a component/unit deteriorates out with time); or negative (in which case, time has a beneficiary effect on the residual lifetime). These
notions of ageing are captured through the well-known monotonic ageing families like Increasing Failure Rate (IFR), Increasing Failure Rate on the Average (IFRA), New Better than Used (NBU), New Better than Used in Expectation (NBUE), Harmonic New Better than Used in Expectation (HNBUE), Decreasing Mean Residual Life (DMRL) and Decreasing Mean Time to Failure (DMTTF)  and their duals. Other important (but perhaps less popular) ageing classes are New Better than Used in Convex ordering (NBUC), new better than renewal used in moment generating function (NBRUmgf),  New Better than Renewal Used in Laplace Transform Order (NBRULC).) However, in practice, one often encounters biological and mechanical systems that exhibit non-monotonic behavior of the mean residual life function or the hazard rate function. Such non-monotonic ageing has often been modelled, for example, using the Bathtub Failure Rate (BFR) function. Another approach using the mean residual life function gives rise to the Increasing then Decreasing Mean Residual Life (IDMRL) class which was introduced by Guess et al. [1]. Later, Mitra and Basu [2] introduced the class of `New Worse then Better than Used in Expectation' (NWBUE) life distributions to model non-monotonic ageing. More recently, Izadi et al. [3] presented the class of Increasing then Decreasing Mean Time To Failure (IDMTTF) to model non-monotonic ageing.

While modelling real life data, the analyst has to check whether ageing has any effect on the characteristic under study. This can be done by checking whether the phenomenon can be modelled by the exponential distribution which implies no ageing. There are many tests for exponentiality, depending on the characterization used. Mention may be made of Xiong et al. [4]  (stemming from the extropy of record values); Rajesh and Sajily [5] (based on the relative extropy); Cuparić et al. [6] (using the notion of the V-empirical Laplace transforms). Cuparić and Milošević [7] are concerned with a new characterization-based of exponentiality tests for randomly censored data. 
 Ossai et al. [8] provide a useful review of tests for exponentiality; while 	Ebner [9] comments on a test for exponentiality based on the MRL function.

Tests for exponentiality against specific monotonic and non-monotonic ageing notions abound in the literature. 
Proschan and Pyke [10] were perhaps the first one to suggest a test for exponentiality against the notion of  IFR alternatives by using the ranks of the normalized spacings between the ordered observations. Later Anis [11] and the references therein also discussed testing exponentiality against IFR alternatives.  IFRA alternatives have been considered by  Ahmad and Mugdadi [12] and the references therein. Testing against NBU alternatives had been proposed initially by Hollander and Proschan [13].  Testing against NBUE alternatives, initially studied by Hollander and Proschan [14], has also been considered by   Anis and Basu [15] as well as Sreelakshmi et al. [16] and the references therein. For testing the notion of exponentiality against HNBUE alternatives (which was first considered by Klefsjö [17,18]), see Bakr et al. [19] and the references therein. A good summary of many of these testing procedures is available in Lai [20]. 

Testing exponentiality against the  recent but somewhat less popular ageing classes has also been extensively explored in the literature. Testing exponentiality against NBUC alternatives has been studied by Mohamed [21]. Hassan et al. [22] and Bakr et al. [23] are concerned with testing exponentilaity against the NBRUmgf Class. Balakrishnan et al. [24] obtained an exact test for testing exponentiality against renewal increasing mean residual life class. Etman et al. [25] explored the problem of testing exponentiality against the   NBRULC class.  In a recent work, Etman et al. [26] obtained a  test for exponentiality against the EBUCL of alternatives. Ghosh and Majumder [27] consider testing exponentiality against DMTTF alternatives using a moment inequality  while Das and Ghosh [28]  present a test for the same problem  based on quantiles. Das and Ghosh [29]  have proposed a class  of tests against DMTTF alternatives. 	Bhattacharyya et al. [30] present a family of tests for the same problem for censored data.

We shall now turn our attention to the problem of testing exponentiality against non-monotonic ageing concepts. Tests for exponentiality against BFR distributions are plenty in the literature; see, for example, Saha and Anis [31] and the references therein. Bhattacharyya et al. [32] extended the BFR notion to the class of bathtub failure rate on the average (BFRA) and considered its testing problem. The problem of testing exponentiality against IDMRL alternatives has been explored by Na and Lee [33] and the references therein.  This class has been extended by Sepehrifar [34] to the recently introduced  Starshaped Mean Equilibrium Life (SMEL) class. 
The problem of testing exponentiality against NWBUE family of life distribution was first discussed by Anis and Mitra [35]. Later Ghosh and Mitra [36] derived the small sample distribution of generalized Hollander-Proschan type statistics for the same testing problem.  Anis [37] offers a survey of tests of exponentiality against non-monotonic ageing classes available till that point of time.

In the hierarchy of monotonic ageing classes, the NBUE class occupies a
special place. It was discussed by Esary et al. [38]. In fact, it was
referred to as \emph{ the net decreasing mean residual life time class} by
Bryson and Siddiqui [39]. Marshall and Proschan [40]  in their
Theorem 2.2 established that the mean waiting time between
consecutive in service system failures is smaller  when a ``failure
replacement policy'' is used compared to an ``age replacement policy'', if and only
if the system life distribution $ F $  has the new better than used in
expectation property (i.e. $ F $ is NBUE). Note that this mean waiting
time is the same under both these policies if the system life can be modelled by the exponential distribution.

If the expected waiting time between the consecutive failures is a principal criterion in choosing whether to endorse an age replacement
policy over the unplanned replacement policy for a given system, then
one possible method to resolve would be to test the above hypothesis for the
life distribution of the given system on the basis of the life data of
the same. Rejection of $ H_{0} $  would indicate  as supporting the adoption of
age replacement plan.

Moment bounds have been derived for many of these non-parametric
families of distributions. Basu and Simons [41] provide lower bounds
for the   $ r-\mathrm{th} $ moment of an IFR distribution. Sengupta and
Deshpande [42] provide moment bounds for a distribution when it is ordered with respect to a known distribution.

Analogy can be drawn between monotonic and non-monotonic ageing classes. The BFR class can be viewed as a non-monotonic version of DFR-IFR classes. Again the IDMRL class is the corresponding non-monotonic version of the IMRL-DMRL classes. Furthermore, the IDMTTF class may be considered as the non-monotonic anlogue of the DMTTF class.  Similarly, the NWBUE class is the non-monotonic genre of the NWUE-NBUE classes. Thus, it should be noted that an NBUE distribution is nothing but an NWBUE distribution with change point $x_{0}=0.$ Results on shock model carry over in an analogous way from the NBUE class to  the NWBUE class of distributions. Similarly, weak convergence is observed within both the classes (as we shall show).

This paper discusses the NBUE class and its non-monotone version, the NWBUE class; and is arranged as follows.  We discuss the NBUE class in Section 2. Specifically, we establish some moment inequalities for this class and prove its weak convergence. Moving forward in Section 3, we focus on the NWBUE class. We 
construct a counterexample to show (as expected) that neither the NBUE moment bounds hold for the NWBUE distributions nor the moment bounds for the NWBUE class of distributions hold for the NBUE class. We also examine  the relationship between the IDMRL and NWBUE classes and in the process correct a long standing fallacy. An erroneous example of NWBUE class of life distribution given in Mitra and Basu [2] is  pointed out and subsequently corrected. Section 4 concludes the paper.

\section{Results on the NBUE Class }

We shall, in this section, state and prove the pertinent results for the
NBUE family of life distributions. We shall first derive some moment inequalities and associated results. Next we shall investigate weak convergence within the NBUE class of life distributions. However, we begin with the definition of the NBUE class for completeness.

\textbf{Definition 1:} Let $ X $ be a non-negative random variable with distribution function $ F $ and finite mean $ \mu. $ $ X $ is said to be New Better than Used in Expectation (NBUE) if 
$$ \int_{0}^{\infty}\overline{F}\left( x|t\right)  dx\leq \mu, \mbox{~~~~~} \forall t\geqslant 0. $$

\subsection{Moment Inequalities and Related Results}
The most important result is given in Theorem  \ref{Th2} below. However, first we shall note the following useful lemma (due to Marshall and Proschan [40]) for completeness:

\begin{lemma} (Marshall and Proschan [40]) \label{L1}
If $ F $ is NBUE with finite mean $ \mu, $ then $ \int_{x}^{\infty} \bar{F}(t) dt\leqslant \mu e^{-x/\mu}, x\geqslant 0, $ where $ \bar{F}=1-F. $
\end{lemma}

We shall now state a well-known general theorem and furnish its proof.
\begin{theorem} \label{Th1}
If $ F $ is NBUE with finite mean $ \mu $ and $ \phi $ is an non-negative non-deceasing function on $ \left[ 0, \infty\right) , $ then 
$$ \int_{0}^{\infty}\phi\left( y\right) \bar{F}\left( y\right) dy \leqslant \int_{0}^{\infty}\phi\left( y\right) e^{-y/\mu} dy.   $$
\end{theorem}
\begin{proof}
Define the "\emph{first derived distribution}" $ f_{\left( 1\right)}$ of $ F $ by $f_{\left( 1\right)}\left( x\right) = \overline{F}\left( x\right) /\mu, \mbox{~~~} x> 0.  $
We have
\begin{eqnarray*}
\int_{0}^{\infty}\phi\left( y\right) \bar{F}\left( y\right) dy & = & \mu \int_{0}^{\infty}\phi\left( y\right) f_{\left( 1\right) }\left( y\right) dy\\
& = & \mu \int_{0}^{\infty}\phi\left( y\right) dF_{\left( 1\right) }\left( y\right) \\
& = & \mu E\left[ \phi \left( Z\right) \right], \mbox{where } Z = F_{\left( 1\right) } \mbox{ in distribution}\\
& \leqslant & \mu \int_{0}^{\infty} P\left[ \phi\left( Z\right) > x\right] dx\\
& = & \mu \int_{0}^{\infty} P\left[  Z > \phi^{-1}\left( x\right) \right] dx, \mbox{where $ \phi^{-1} $ is a suitably defined inverse function of $ \phi $}\\
& \leqslant & \mu \int_{0}^{\infty} P\left[  U > \phi^{-1}\left( x\right) \right] dx, \mbox{where $ P\left[ U>t\right] = e^{-t/\mu}, $ using Lemma \ref{L1}}\\
& = & \mu \int_{0}^{\infty} P\left[ \phi\left( U\right) > x\right] dx\\
& = & \mu E \left[\phi\left( U\right)\right] \\
& = & \mu \int_{0}^{\infty}\phi\left( y\right) \frac{1}{\mu}e^{-y/\mu}dy\\
&=&  \int_{0}^{\infty}\phi\left(y\right) e^{-y/\mu}dy.
\end{eqnarray*}
\end{proof}
We now come to the main result given below as Theorem \ref{Th2}.

\begin{theorem}\label{Th2}
If $ F $ is NBUE with finite mean $ \mu,$ then 
\[\mu_{r}=E\left( X^{r}\right) \left\lbrace 
\begin{array}{ll}
\leqslant \Gamma\left( r+1\right) \mu^{r} & \mathrm{if} \mbox{~~} r\geqslant 1,\\
\geqslant \Gamma\left( r+1\right) \mu^{r} & 0<r<1.
\end{array}
\right. 
\]
\end{theorem}

\begin{proof}
The proof is  established easily by appealing to Theorem \ref{Th1}, with $ \phi\left( x\right) =x^{r-1}. $
\end{proof}

\subsection{Weak Convergence within the NBUE Class }
Basu and Simons [41] had shown the weak convergence of the IFR class. Subsequently, Basu and Bhattcharjee [43] established a similar conclusion for the \mbox{HNBUE} class. Since the HNBUE class is strictly larger than the NBUE class, it follows that if $ F_{n}, n=1, 2, 3,\cdots  $ be a sequence of NBUE life distributions converging to $F,$ then $F$ is HNBUE. However, here, we establish a more powerful result, namely that the limiting distribution $ F $ is indeed NBUE.  We begin with the following theorem.
\begin{theorem}\label{Th3}
Let $ F_{n}, n=1, 2, 3,\cdots  $ be a sequence of NBUE life distributions with mean $ \mu_{n}.$ Assume that
\begin{itemize}
\item[(a)]
 $F_{n} \rightarrow F$ \mbox{ in law, where }  $F$  \mbox{is a continuous distribution function with finite mean }  $\mu;$  \label{a1}

 and
\item[(b)] the sequence $ \left\lbrace \mu_{n}\right\rbrace  $ is bounded.
\end{itemize}
Then $ F $ is NBUE. Furthermore, 
\begin{equation} 
\lim_{ n\rightarrow \infty} \int_{0}^{\infty}x^{r}dF_{n}\left( x\right) = \int_{0}^{\infty}x^{r}dF\left( x\right) \label{Eq1}
\end{equation}
for every $ r>0. $
\end{theorem}
\begin{proof}
Let $ \mu_{n}\leqslant M, \forall n\geqslant 1. $ By appealing to Fatou's lemma together with the assumptions of Theorem \ref{a1}, we have  
\begin{equation}
\mu=\int_{0}^{\infty}\overline{F}\left( x\right)  dx\leqslant M < \infty. \label{Eq2}
\end{equation}

We shall first show that $ \mu_{n}\rightarrow \mu $ as $ n\rightarrow \infty . $\\
For any constant $ C >0, $ to be chosen suitably, we have
\begin{eqnarray*}
\int_{0}^{\infty}\left\lbrace  \overline{F}_{n}\left( x\right)  - F_{n}\left( x\right)\right\rbrace  dx & = & \int_{0}^{C}\left\lbrace  \overline{F}_{n}\left( x\right)  - F_{n}\left( x\right)\right\rbrace  dx + \int_{C}^{\infty}\overline{F}_{n}\left( x\right)dx - \int_{C}^{\infty}\overline{F}\left( x\right)dx\\
& = & I_{1n}+I_{2n}-I_{3},
\end{eqnarray*}
where
$$ I_{1n} = \int_{0}^{C}\left\lbrace  \overline{F}_{n}\left( x\right) - F_{n}\left( x\right)\right\rbrace  dx;   $$
$$I_{2n}= \int_{C}^{\infty}\overline{F}_{n}\left( x\right)dx,  $$
and 
$$I_{3}= \int_{C}^{\infty}\overline{F}\left( x\right)dx. $$
First observe that under the given conditions, by appealing to the Dominated Convergence Theorem, we get $I_{1n}\rightarrow 0  $ as $ n\rightarrow\infty. $\\
Secondly, we have,  
\begin{eqnarray}
I_{2n}=\int_{C}^{\infty}\overline{F}_{n}\left( x\right)dx & \leqslant & \mu_{n}\exp \left( -C/\mu_{n}\right) (\mbox{by  Lemma 2.1)}\\
& \leqslant & M\exp \left( -C/M\right) 
\end{eqnarray}
which can be made as close to  to zero as we want by choosing $ C $ sufficiently large.\\
Finally, because of (\ref{Eq2}), $ I_{3}= \int_{C}^{\infty}\overline{F}\left( x\right)dx $ can be made as close to zero as we want by choosing $ C $ sufficiently large.
Thus, $ \int_{0}^{\infty}\left\lbrace  \overline{F}_{n}\left( x\right)  - F_{n}\left( x\right)\right\rbrace  dx \rightarrow 0,  $ as $ n\rightarrow \infty. $ Hence, in other words,
\begin{equation}
 \mu_{n}\rightarrow \mu  \mbox{ as }  n\rightarrow \infty. \label{a3}
\end{equation}

Thus,
\[\frac{1}{\overline{F}_{n}\left( x\right) }\int_{x}^{\infty}\overline{F}_{n}\left( t\right)dt\leqslant \mu_{n},\]
since $ F_{n} $ is NBUE. Hence, by taking limits as $ n\rightarrow \infty, $ and using the Dominated Convergence Theorem, we get 
\[\frac{1}{\overline{F}\left( x\right) }\int_{x}^{\infty}\overline{F}\left( t\right)dt\leqslant \mu.\]
Thus, $ F $ is NBUE. This settles the proof of the first part of the theorem.\\

Now to prove (\ref{Eq1}), we observe that using (\ref{a3}), the result holds for $ r=1.$ Hence, (see e.g. Chung [44], p. 88), (\ref{Eq1}) holds for $ r<1 $ as well.\\
For $ r>1, $ by Theorem \ref{Th2}, it follows that

\[\mu_{n;r}=E\left( X_{n}^{r}\right) \leqslant \mu_{n}^{r}\Gamma \left( r+1\right). \]
Hence,
\begin{equation}
 \sup_{n} E\left( X_{n}^{r}\right)<\infty, \forall r>1. 
\end{equation}  
Thus, $ \left\lbrace X_{n}^{r}\right\rbrace  $ is uniformly integrable for each $ r>1. $ This proves (\ref{Eq1}).\\

\textbf{Remark:} In passing it should be mentioned that there is another elegant way to prove (\ref{Eq1}). Note that since $ X_{n} $ is NBUE, it admits finite moments of all orders (by Theorem \ref{Th2}). Hence, by appealing to Corollary to Theorem 25.12 of Billingsley [45], (\ref{Eq1}) follows.
\end{proof}

\begin{lemma}\label{L2}
A distribution function which is NBUE is uniquely determined by its moment sequence.
\end{lemma}
\begin{proof}
Let $ F $ be NBUE with mean $ \mu. $ In view of Theorem \ref{Th2} and a result from  Lo\'{e}ve [46], p 217, it follows that an NBUE distribution is uniquely determined by its moment sequence.
\end{proof}

We now give the converse of Theorem \ref{Th3}.
\begin{theorem}\label{TH4}
Let $ F_{n}, n=1, 2, 3,  \cdots $ be a sequence of NBUE distributions, such that for each integer $ r>0, $
\begin{equation} 
\lim_{ n\rightarrow \infty} \int_{0}^{\infty}x^{r}dF_{n}\left( x\right) = \int_{0}^{\infty}x^{r}dF\left( x\right) \label{Eq1a}
\end{equation}
for some life distribution $ F. $ Then $ F_{n} \rightarrow F$ in law, and $ F $ is NBUE.
\end{theorem}
\begin{proof}
By Lemma \ref{L2}, $ F $ is uniquely determined by its moment sequence. Hence, using (\ref{Eq1a}), the limiting distribution of every weakly convergent subsequence of $ \left\lbrace F_{n}\right\rbrace  $ is necessarily $ F. $ Therefore, using a conventional  argument based on tightness and relative compactness together with the fact that an NBUE life distribution is uniquely determined by its moment sequence, completes the proof of the theorem. 
\end{proof}

\section{Results on the NWBUE Class }
We shall now turn our attention to the NWBUE family of life distributions. We repeat the following definition for completeness.

\textbf{Definition 2:} Let $ X $ be a non-negative random variable with cumulative distribution function $ F $ and finite mean $ \mu. $ $ X $ is said to be New Worse then Better than Used in Expectation (NWBUE) if there exists a point $ t_{0}\geqslant 0 $ such that the mean residual life at age $ t, $  (denoted by $ e_{F}\left( t\right)  $)
\[e_{F}\left( t\right)  \left\lbrace 
\begin{array}{ll}
\geqslant e_{F}\left( 0\right) & for \mbox{~~} t < t_{0} ,\\
\leqslant e_{F}\left( 0\right)  & for \mbox{~~} t \geqslant t_{0}.
\end{array}
\right. 
\]Observe that the NWBUE class can be viewed as the non-monotonic analogue of the NBUE class.   First, we shall show that the moment bounds for the NBUE class do not hold for the NWBUE class and vice-versa.  Secondly, we shall modify and correct a result involving IDMRL and NWBUE classes of life distributions and finally point out and correct an error in Mitra and Basu [2].
\subsection{Moment Bounds for NWBUE Distributions } 
 It is natural to expect that the moment bounds for NWBUE distributions involve the change point $ x_{0}. $ These bounds have been obtained by Mitra and Basu [2]. For completeness, we reproduce their results below.

\begin{theorem}(Mitra \& Basu [2])
If $ F $ is $ NWBUE\left( x_{0}\right) , x_{0} < \infty,\)$ 
$
with
finite mean $ \mu $  then

\begin{itemize}

\item[$\mathrm{(a)}$] \[\mu_{r} = E\left( X^{r} \right)\leq \left( \geqslant\right)  \text{\ r}e^{\frac{x_{0}}{\mu}}\int_{x_{0}}^{\infty}{x^{r - 1}e^{- \frac{x}{\mu}}\text{dx}\text{\ \ \ \ \ \ for\ }r \geqslant 1 \left( < 1\right);}\]

\item[$\mathrm{(b)}$]\[\mu_{r} \leq x_{0}^{r} + \mu^{r}\Gamma\left( r + 1 \right)\sum_{j = 0}^{r - 1}\frac{\left( \frac{x_{0}}{\mu} \right)^{r}}{r!}\text{\ \ for\ all\ integers\ }r \geq 1;\]

\item[$\mathrm{(c)}$]\[\mu_{r} \leq \mu^{r}\Gamma\left( r + 1 \right)e^{\frac{x_{0}}{\mu}}\text{\ \ for\ all\ \ }r \geq 1.\]
\end{itemize}

\end{theorem}

We shall give  counter-examples to show that  the moment bounds for NWBUE distributions do not hold for NBUE distributions and vice-versa. Towards this end, we construct the
following example.

\textbf{Example 3.1:}  Consider an MRL function defined as follows :-

\[e_{F}\left( x \right) = \left\{ \begin{matrix}
5 + x\ \mbox{~~~} 0 \leq x \leq 1 \\
\frac{55 - x}{9} \mbox{~~~}1 \leq x \leq 28 \\
3 \mbox{~~~}x\geq 28. \\
\end{matrix} \right.\ \]

Here the mean \(\mu = e_{F}\left( 0 \right) = 5.\) It is easy to see
that the above distribution is NWBUE with the corresponding change point \(x_{0} = 10.\)
Recall that the MRL function \(e_{F}\left( x \right)\) uniquely
determines the distribution function \(F\left( x \right);\) see for
example Berilant et al. [47]. The relationship between the survival
function \(\overline{F}\left( x \right)\) and the MRL function
\(e_{F}\left( x \right)\) is given by

\[\overline{F}\left( x \right) = \frac{e_{F}\left( 0 \right)}{e_{F}\left( x \right)}\exp\left\lbrack - \int_{0}^{x}\frac{1}{e_{F}\left( u \right)}du \right\rbrack.\]

Hence by using the above mentioned inverse relationship, it follows that

\[\overline{F}\left( x \right) = \left\{ \begin{matrix}
\frac{25}{\left( 5 + x \right)^{2}} \mbox{~~~~~~~~~~~~~~~~~~~}0\leq x \leq 1 \\
\frac{75}{2 \times 54^{9}}\left( 55 - x \right)^{8}\mbox{~~~~~~~~~~~}1 \leq x \leq 28 \\
\frac{25}{9 \times 2^{10}}\exp\left\lbrack \frac{1}{3}\left( 28 - x \right) \right\rbrack \mbox{~~~}x \geq 28. \\
\end{matrix} \right.\ \]

It is easy to examine that $ F $  is indeed a distribution function. Then
routine calculation gives \(E_{F}\left( X^{2} \right) = 54.1210;\) and
all the three bounds given by Mitra and Basu [2]  are satisfied.
But \(\Gamma\left( 2 + 1 \right)\mu^{2} = 2 \times 5^{2} = 50.\) Thus
clearly

\[D\left( t \right) = \mu^{t}\Gamma\left( t + 1 \right) - E\left( X^{t} \right) < 0,\ \text{for\ }t = 2 > 1,\]

thereby showing  that the moment bounds
for NWBUE distributions do not hold good for NBUE distributions.\\

\textbf{Example 3.2:} Let $ X $ be a Weibull random variable with scale parameter unity and shape parameter 2. By virtue of this choice of the shape parameter, $ X $ is IFR and hence NBUE.  By Theorem \ref{Th2},  for $ r\geqslant 1 $ we have, $ \mu_{r}= E\left( X^{r}\right)\leqslant \mu^{r}\Gamma\left( r+1)\right) >0.$ In particular, for $ r=2,$ we have $ \mu_{2}=1,$ while the corresponding NBUE moment bound is $ \frac{\pi}{2}=1.57.$ The corresponding bound for NWBUE distributions is $\frac{{\pi}^{3/2}}{4}=1.392.$ Thus, this  random variable belonging to the NBUE class does not satisfy the moment bound for NWBUE distributions. \

These two examples show that the bound on moments do \emph{not} carry over from NBUE class to the corresponding non-monotonic NWBUE class.

\subsection{Relationship with IDMRL Class}
Mitra and Basu [2] have proved in their work  that under appropriate regularity conditions, all BFR distributions are NWBUE distributions. The connection with IDMRL class is \emph{incorrectly} investigated in their Proposition 2.3, which is reproduced below for completeness. \\

\textbf{Proposition 2.3 }(Mitra and Basu [2]). 
If $ F $ is $ IDMRL\left( t_{0}\right) , $ then $ F $ is $NWBUE\left( t_{0}^{\star}\right)  $ with $t_{0}^{\star}\geqslant t_{0}. $

Mitra and Basu [2] omitted its proof ``since it is trivial''. Unfortunately, this notion of NWBUE class containing \emph{all} IDMRL distributions is quiet prevalent in the statistical literature. See, for example,  Lai and Xie [48] and Anis [49].  However, we give below an example of an IDMRL distribution which is \emph{not} NWBUE.
 
\textbf{Example 3.3:} Let $ \overline{G}$ be a survival function defined as follows:-
\[\overline{G}\left( x\right) =\left\lbrace 
  \begin{array}{l l}
  \frac{4}{\left( x+2\right) ^{2}} & 0\leqslant x \leqslant 2\\
  \frac{10-x}{32} & 2\leqslant x \leqslant 4\\
  \frac{3}{16}\exp \left( \frac{4-x}{3}\right)  & x\geqslant 4.
    \end{array}
    \right.
\]
Its associated MRL function is given by 

\[e_{G}\left( x \right) =\left\lbrace 
  \begin{array}{l l}
  x+2  & 0\leqslant x \leqslant 2\\
  \frac{10-x}{2} & 2\leqslant x \leqslant 4\\
  3  & x\geqslant 4.
    \end{array}
    \right.
\]

It is easy to show that $ G $ is IDMRL with change point $ \tau_{0}=2. $ The mean of $ G $ is clearly $ \mu=2; $ thus, the MRL function is always above the mean and \emph{never} crosses it. Thus $ G $ is IDMRL; but not NWBUE. We correct the Proposition 2.3 of Mitra and Basu [2] by presenting the following theorem. 

\begin{theorem}
Let $ F $ be an IDMRL (DIMRL) distribution with change point $ \tau_{0}, $ $ F^{-1}\left( 0\right) \leqslant \tau_{0} \leqslant  F^{-1}\left( 1\right). $ If there exists an $ x^{\star} \in \left(\tau_{0}, F^{-1}\left( 1\right)\right) $ such that $e_{F}\left( x^{\star} \right) =\mu,$ then $ F $ is NWBUE (NBWUE) with change point $x^{\star}.$ If there does not exist such an $x^{\star},$ then $ F $ is NWUE (NBUE).
\end{theorem}
\begin{proof}
Assume that $ F $ is IDMRL with change point  $ \tau_{0}.$ Since $e_{F}\left( 0 \right) =\mu$ and $e_{F}\left( x\right)$ is increasing on $\left( F^{-1}\left( 0\right), \tau_{0}\right),$ it is obvious that $ e_{F}\left(x\right) \geqslant \mu $ for every $ x \in \left[ F^{-1}\left( 0 \right),\tau_{0}\right).$ Since $ e_{F}\left( x\right)$ is decreasing on $ \left( \tau_{0}, F^{-1}\left( 1\right) \right),  $ the function $e_{F}\left( x\right)- \mu$ changes sign at most once from `+' to `-'on $ \left(\tau_{0},F^{-1}\left( 1\right) \right).$  If it does, then $ F $ is NWBUE; if it does not then $ F $ is NWUE.\\
The proof of the dual is similar.
\end{proof}

\subsection{A Correct Example}
We shall close this section  by presenting a correct example of an NWBUE distribution, thereby rectifying the example presented in Mitra and Basu [2].

The first example of an NWBUE distribution given in Mitra and Basu [2] (their Example 1.4) is \emph{incorrect}, since $ F $ is not a cumulative distribution function. Observe that their $ \overline{F} $ is increasing in $ x $ for $ x \in \left[ 1, \sqrt{3}\right); $ but a survival function (which $ \overline{F} $ is supposed to represent), cannot be an increasing function. Hence, their $ e_{F}\left( x\right)  $ does not represent an MRL function. There are certain necessary conditions to be fulfilled for a function to be a valid MRL function. These necessary conditions may be found in Gupta and Kiramni [50], Bhattacharjee [51] and Guess and Proschan [52].  It is easy to check that the corresponding $ e_{F}\left( x\right)  $ as given in Mitra and Basu [2] does not satisfy the necessary conditions to represent an MRL function of a random variable non-degenerate at 0 as given in Gupta and Kiramni [50], thereby confirming that it is not a valid MRL function. We, therefore, correct their example as follows.

\textbf{Example 3.4:} Consider the life distribution whose survival function is given by

 \[
 \bar{F}\left( x\right)=\left\lbrace
 \begin{array}{ll}
 \frac{4}{\left( 2+x\right) ^{2}} & 0\leqslant x<1 \\
 \frac{2}{27}\left( 7-x\right)  & 1\leqslant x <3\\
 \frac{2x\left( 3+x\right) ^{2}}{729}\exp \left( 3-x\right) & 3 \leqslant x <\infty .
\end{array}  
\right. 
 \]
The corresponding MRL function is given by 

\[ e_{F}\left( x\right)=\left\lbrace 
\begin{array}{ll}
2+x & 0\leqslant x<1 \\
\frac{7-x}{2} & 1\leqslant x <3\\
1+ \frac{3}{x} & 3 \leqslant x <\infty .
\end{array}
\right.
 \]
It is easy to check that $\bar{F}\left( x\right)  $ is indeed a survival function; and $e_{F}\left( x\right)  $ is an MRL function satisfying all the conditions mentioned in Gupta and Kiramni [50]. Moreover, it is easy to see that $ F $ is NWBUE with change point $ x_{0}=3. $

\section{Conclusion}
In the context of the NBUE family, we prove its weak convergence. The NWBUE class of distributions is the non-monotonic analogue of the
NBUE class of distributions. While results relating to shock model
theory carry over from the monotonic class to the non-monotonic class in
an appropriate manner, it is shown by means of a counter example that
the same does not hold true for moment bounds. We also correct a long standing mistaken notion regarding the relationship between the IDMRL and NWBUE classes.


\textbf{References}

\begin{enumerate}
\def\labelenumi{\arabic{enumi}.}

\item
  Guess, F., Hollander, M. \& Proschan, F. (1986). Testing exponentiality
  versus a trend change in mean residual life. \emph{Annals of
  Statistics}. \textbf{14}; 1388-98.
  
\item
  Mitra, M. \& Basu, S. K. (1994). On a non-parametric family of life
  distributions and its dual. \emph{Journal of Statistical Planning and
  Inference.}\textbf{39}; 385-397. 
  
\item
Izadi, M., Sharafi, M. and Khaledi, B.E., (2018). New nonparametric classes of distributions in terms of mean time to failure in age replacement. \textit{Journal of Applied Probability}, \textbf{55}(4), 1238-1248.  

\item
Xiong, P., Zhuang, W., \& Qiu, G. (2022). Testing exponentiality based on the extropy of record values. Journal of Applied Statistics, 49(4), 782-802.

\item
Rajesh, G., Sajily, V.S. Testing exponentiality based on relative extropy. Stat Papers 66, 108 (2025). https://doi.org/10.1007/s00362-025-01726-6

\item
Cuparić, M., Milošević, B., \& Obradović, M. (2022). New consistent exponentiality tests based on V-empirical Laplace transforms with comparison of efficiencies. Revista de la Real Academia de Ciencias Exactas, Físicas y Naturales. Serie A. Matemáticas, 116(1), 42.

\item
Cuparić, M., \& Milošević, B. (2022). New characterization-based exponentiality tests for randomly censored data. Test, 31(2), 461-487.

\item
  Ossai, E. O.; Madukaife, M. S. and Oladugba, A. V. (2022). A review of
  tests for exponentiality with Monte Carlo comparisons. \emph{Journal
  of Applied Statistics.}\textbf{49}(5), 1277-1304. doi.org/10.1080/02664763.2020.1854202

\item  
Ebner, B. (2023). The test of exponentiality based on the mean residual life function revisited. Journal of Nonparametric Statistics, 35(3), 601-621.

\item  
  Proschan, F.,\& Pyke, R. (1967). Tests for monotone failure rate. Proceedings of the Fifth Berkeley Symposium on Probability and Mathematical Statistics 3, 293–312.
  
\item
  Anis, M. Z. (2013). A family of tests for exponentiality against IFR
  alternatives. \emph{Journal of Statistical Planning and Inference}.
  \textbf{143}; 1409-1415.
  
\item
  Ahmad, I. A. \& Mugdadi, A. R. (2004). Further moment inequalities of
  life distributions with hypotheses testing applications: the IFRA,
  NBUC and DMRL classes. \emph{Journal of Statistical Planning and
  Inference}. \textbf{120}; 1-12.
  
\item
  Hollander, M. \& Proschan, F. (1972). Testing whether new is better
  than used. \emph{Annals of Mathematical Statistics.}\textbf{43};
  1136-1146.
\item
  Hollander, M. \& Proschan, F. (1975). Tests for the mean residual life.   \emph{Biometrika}. \textbf{62}; 583-593.

\item
  Anis, M. Z. \& Basu, K. (2011). The exact null distribution of the generalized Hollander–Proschan type test for NBUE alternatives.  \emph{Statistics \& Probability Letters.}   \textbf{81} (11); 1733-1737.

\item  
  Sreelakshmi, N., Kattumannil, S. K., \& Asha, G. (2018). Quantile based tests for exponentiality against DMRQ and NBUE alternatives. Journal of the Korean Statistical Society, 47(2), 185-200.
  
\item
  Klefsjö, B. (1983a).Some tests against ageing based on the total time
  on test transform. \emph{Communications in Statistics --Theory and
  Methods}. \textbf{12}; 907-927.
\item
  Klefsjö, B. (1983b). Testing exponentiality against HNBUE.
  \emph{Scandinavain Journal of Statistics.} \textbf{10}; 67-75.
  
\item  
Bakr, M. E., Kibria, B. G., \& Gadallah, A. M. (2023). A new non-parametric hypothesis testing with reliability analysis applications to model some real data. Journal of Radiation Research and Applied Sciences, 16(4), 100724.  

\item
  Lai, C. D. (1994). Tests of univariate and bivariate stochastic
  ageing. \emph{IEEE Transactions on Reliability.}\textbf{43}; 233--240. 
\item  
Mohamed, N. M. (2024). Testing Exponentiality Against NBUC Class with Some Applications in Health and Engineering Fields. Frontiers in Scientific Research and Technology, 8(1).

\item
Hassan, N. A., Said, M. M., Attwa, R. A. E. W., \& Radwan, T. (2024). Applying the Laplace Transform Procedure, Testing Exponentiality against the NBRUmgf Class. Mathematics, 12(13), 2045. https://doi.org/10.3390/math12132045

\item 
Bakr, M. E., El-Atfy, E. S., Balogun, O. S., Tashkandy, Y. A., \& Gadallah, A. M. (2024). Statistical Insights: Analyzing Shock Models, Reliability Operations and Testing Exponentiality for NBRUmgf Class of Life Distributions. Maintenance \& Reliability/Eksploatacja i Niezawodność, 26(2).

\item
Balakrishnan, N., Mathew, D. C., \& Kattumannil, S. K. (2022). An exact test for exponentiality against renewal increasing mean residual life class. Statistics, 56(1), 164-181.

\item
Etman, W. B. H., EL-Sagheer, R. M., Abu-Youssef, S. E., \& Sadek, A. (2022). On some characterizations to NBRULC class with hypotheses testing application. Appl. Math. Inf. Sci, 16, 139-148. 

\item
Etman, W. B., Abouelenein, M. F., Eliwa, M. S., El-Morshedy, M., Roushdy, N., \& EL-Sagheer, R. M. (2025). A Novel U-Statistic Test for Exponentiality Against EBUCL Reliability Class and Applied to Complete and Censored Data Across Various Risk Profiles. European Journal of Pure and Applied Mathematics, 18(3), 6376-6376. 

\item
Ghosh, S., \& Majumder, P. (2022). A moment inequality for decreasing mean time to failure distributions with hypothesis testing application. Journal of Statistical Computation and Simulation, 92(14), 2875–2890.\\ https://doi.org/10.1080/00949655.2022.2051172

\item
Das, K. and  Ghosh, S.(2025)A Quantile-Based Approach for Testing Exponentiality Against DMTTF Alternatives. IEEE Transactions on Reliability, doi: 10.1109/TR.2025.3560678.

\item
Das, K., \& Ghosh, S. (2025). A class of nonparametric tests for DMTTF alternatives based on moment inequality. Statistical Papers, 66(2), 47\\
https://doi.org/10.1007/s00362-025-01670-5

\item
Bhattacharyya, D., Khan, R. A., \& Mitra, M. (2023). A new family of tests for DMTTF alternatives under complete and censored samples. Communications in Statistics-Simulation and Computation, 52(10), 4603-4620.DMTTF

\item
Saha, A., \& Anis, M. Z. (2023). A family of tests for trend change in failure rate function with right censored data. \textit{Journal of Statistical Computation and Simulation}. \textbf{94}(6), 1191–1203. https://doi.org/10.1080/00949655.2023.2282740  

\item
Bhattacharyya, D., Ghosh, S., \& Mitra, M. (2022). On a non-monotonic ageing class based on the failure rate average. Communications in Statistics-Theory and Methods, 51(14), 4807-4826.

\item
  Na, M. H. \& Lee, S. (2003). A family of IDMRL tests with unknown
  turning point. \emph{Statistics}. \textbf{37}(5); 457-462. 

\item  
Sepehrifar, M. (2025). A U-statistic-based test for exponentiality using starshaped mean equilibrium class of life distributions. Statistical Papers, 66(6), 1-19.

\item
  Anis, M. Z. \& Mitra, M. (2005). A simple test of exponentiality
  against NWBUE family of life distributions\emph{. Applied Stochastic
  Models in Business and Industry}. \textbf{21}; 45-53.
  
\item  
Ghosh, S., \& Mitra, M. (2022). On the exact distribution of generalized Hollander-Proschan type statistics. Communications in Statistics-Simulation and Computation, 51(9), 5051-5067. 

\item
  Anis, M. Z. (2014). Tests of non-monotonic stochastic ageing
  notions in reliability theory. \emph{Statistical Papers.}55;
  691--714.
  
\item  
Esary, J. D., Marshall, A. W., \& Proschan, F. (1970). Some reliability applications of the hazard transform. \textit{SIAM Journal on Applied Mathematics.} \textbf{18}(4); 849-860.  
  
\item
  Bryson, M. C. \& Siddiqui, M. M. (1969). Some criteria for ageing.
  \emph{Journal of the American Statistical Association}.
  \textbf{64}(328); 1472-1483.
  
\item
  Marshall, A.W. \& Proschan, F. (1972). Classes of distributions
  applicable in replacement with renewal theory implications. In: LeCam,
  L., Neyman, J., Scott, E.L. (Eds.), \emph{Proceedings of the Sixth
  Berkeley Symposium of Mathematical Statistics and Probability}, Vol.
  I. University of California Press, Berkeley, CA, 395-415.
  
\item
  Basu, S. K. \& Simons, G. D. (1982). Moment spaces for IFR
  distributions, applications and related material. \emph{Technical
  Report.} University of North Carolina.

\item
  Sengupta, D. \& Deshpande, J . V. (1994). Some results on the relative
  ageing of two life distributions. \emph{Journal of Applied
  Probability.}\textbf{31} (4); 991- 1003.

\item
 Basu, S. K. \& Bhattacharjee, M. C. (1984). On weak convergence within the HNBUE family of life distributions. \textit{Journal of Applied Probability.}  \textbf{21}; 654-660.
 
\item
Chung, K. L. (1968). \textit{A Course on Probability Theory.} Harcourt, Brace and World. New York.

\item
Billingsley, P. (1986). \textit{Probability and Measure,} Second Edition. John Wiley \& Sons, New York.
  
\item
Lo\'{e}ve, M. (1963). \textit{Probability Theory,} 3rd Edition. Van Nostrand. New York.

\item
  Berilant, J.; Broniatowski, M.; Teugels, J. L. \&Vynckier, P. (1995).
  The mean residual life function at great age: Application to tail
  estimation. \emph{Journal of Statistical Planning and Inference}.
  \textbf{45}; 21-48.
  
\item  
  Lai, C. D.  \& Xie, M. (2006). \textit{Stochastic Ageing and Dependence for Reliability.} Springer. New York.

\item
Anis, M. Z. (2012). On some properties of the IDMRL class of life distributions. \textit{Journal of Statistical Planning and Inference}, \textbf{142}(11), 3047-3055.

 \item 
  Gupta, R.C., Kirmani, S.N.U.A. (1998). Residual life function in reliability studies. In: \textit{Frontiers in Reliability}; Eds. A. P. Basu; S. K. Basu \& S. P. Mukherjee; World Scientific, New Jersey, pp. 175–190.
  
 \item 
  Bhattacharjee, M. C. (1982) The class of mean residual lives and some consequences. SIAM J.
 Alg. Disc. Meth. 3, 56-65.
  
\item 
Guess, F. \& Proschan, F. (1988). Mean residual life: theory and applications. In \textit{Handbook of Statistics.} Vol. 7, Ed. P. R. Krishnaiah and C. R. Rao, Elsevier, Amsterdam, pp 215-224.
  
\end{enumerate}

\hrule

\end{document}